\documentclass[a4paper,11pt]{article}
\usepackage[english,russian]{babel}
\usepackage[cp1251]{inputenc}
\usepackage{mathtext}
\usepackage[T2A]{fontenc}
\usepackage{inputenc,amssymb}
\usepackage{epsfig}
\usepackage{caption2}
\usepackage{amsfonts}
\usepackage{amssymb}
\usepackage{amsmath}
\usepackage{geometry}
\usepackage{color}

\newcounter{eqcounter}[section]

\renewcommand{\r}{\mathbb R}
\newcommand{\eqd}{\stackrel{d}{=}}

\title{Generalized negative binomial distributions as mixed geometric laws
and related limit theorems\thanks{Research supported by the Russian
Science Foundation, project 18-11-00155.}}
\author{V. Yu. Korolev\footnote{Faculty of
Computational Mathematics and Cybernetics, Lomonosov Moscow State
University; Institute of Informatics Problems, Federal Research
Center <<Informatics and Control>>, Russian Academy of Sciences;
Hangzhou Dianzi University; victoryukorolev@yandex.ru}, A. I.
Zeifman\footnote{Vologda State University; Institute of Informatics
Problems, Federal Research Center <<Informatics and Control>>,
Russian Academy of Sciences;  Vologda Research Center of the Russian
Academy of Sciences; a$\_$zeifman@mail.ru}}

\date{}

\begin{document}

\maketitle


{\bf Abstract.} In this paper we study a wide and flexible family of
discrete distributions, the so-called generalized negative binomial
(GNB) distributions that are mixed Poisson distributions in which
the mixing laws belong to the class of generalized gamma (GG)
distributions. The latter was introduced by E. W. Stacy as a special
family of lifetime distributions containing gamma, exponential power
and Weibull distributions. These distributions seem to be very
promising in the statistical description of many real phenomena
being very convenient and almost universal models for the
description of statistical regularities in discrete data. Analytic
properties of GNB distributions are studied. A GG distribution is
proved to be a mixed exponential distribution if and only if the
shape and exponent power parameters are no greater than one. The
mixing distribution is written out explicitly as a scale mixture of
strictly stable laws concentrated on the nonnegative halfline. As a
corollary, the representation is obtained for the GNB distribution
as a mixed geometric distribution. The corresponding scheme of
Bernoulli trials with random probability of success is considered.
Within this scheme, a random analog of the Poisson theorem is proved
establishing the convergence of mixed binomial distributions to
mixed Poisson laws. Limit theorems are proved for random sums of
independent random variables in which the number of summands has the
GNB distribution and the summands have both light- and heavy-tailed
distributions. The class of limit laws is wide enough and includes
the so-called generalized variance gamma distributions. Various
representations for the limit laws are obtained in terms of mixtures
of Mittag-Leffler, Linnik or Laplace distributions. Limit theorems
are proved establishing the convergence of the distributions of
statistics constructed from samples with random sizes obeying the
GNB distribution to generalized variance gamma distributions. Some
applications of GNB distributions in meteorology are discussed.

\smallskip


\section{Introduction}

\subsection{Motivation}

In this paper we study a wide and flexible family of discrete
distributions, the so-called generalized negative binomial (GNB)
distributions that are mixed Poisson distributions in which the
mixing laws belong to the class of generalized gamma distributions.
These distributions seem to be very promising in the statistical
description of many real phenomena being very convenient and almost
universal models.

It is necessary to explain why we consider this combination of the
mixed and mixing distributions. First of all, the Poisson mixing
kernel is used for the following reasons. Pure Poisson processes can
be regarded as best models of stationary (time-homogeneous) chaotic
flows of events \cite{BeningKorolev2002}. Recall that the
attractiveness of a Poisson process as a model of homogeneous
discrete stochastic chaos is due to at least two reasons. First,
Poisson processes are point processes characterized by that time
intervals between successive points are independent random variables
with one and the same exponential distribution and, as is well
known, the exponential distribution possesses the maximum
differential entropy among all absolutely continuous distributions
concentrated on the nonnegative half-line with finite expectations,
whereas the entropy is a natural and convenient measure of
uncertainty. Second, the points forming the Poisson process are
uniformly distributed along the time axis in the sense that for any
finite time interval $[t_1,t_2]$, $t_1<t_2$, the conditional joint
distribution of the points of the Poisson process which fall into
the interval $[t_1,t_2]$ under the condition that the number of such
points is fixed and equals, say, $n$, coincides with the joint
distribution of the order statistics constructed from an independent
sample of size $n$ from the uniform distribution on $[t_1,t_2]$
whereas the uniform distribution possesses the maximum differential
entropy among all absolutely continuous distributions concentrated
on finite intervals and very well corresponds to the conventional
impression of an absolutely unpredictable random variable (r.v.)
(see, e. g., \cite{GnedenkoKorolev1996, BeningKorolev2002}). But in
actual practice, as a rule, the parameters of the chaotic stochastic
processes are influenced by poorly predictable ``extrinsic'' factors
which can be regarded as stochastic so that most reasonable
probabilistic models of non-stationary (time-non-homogeneous)
chaotic point processes are {\it doubly stochastic Poisson
processes} also called {\it Cox processes} (see, e. g.,
\cite{Grandell1976, Grandell1997, BeningKorolev2002}). These
processes are defined as Poisson processes with stochastic
intensities. Such processes proved to be adequate models in
insurance \cite{Grandell1976, Grandell1997, BeningKorolev2002},
financial mathematics \cite{KCKZ2015}, physics
\cite{KorolevSkvortsova} and many other fields. Their
one-dimensional distributions are mixed Poisson.

In order to have a flexible model of a mixing distribution which is
``responsible'' for the description of statistical regularities of
the manifestation of ``outer'' stochastic factors we suggest to use
the {\it generalized gamma $($GG$)$ distributions} defined by the
density
$$
g^*(x;r,\alpha,\lambda)=\frac{|\alpha|\lambda^r}{\Gamma(r)}x^{\alpha
r-1}e^{-\lambda x^{\alpha}},\ \ \ \ x\ge0,
$$
with $\alpha\in\mathbb{R}$, $\lambda>0$, $r>0$. The class of GG
distributions was first described as a unitary family in 1962 by E.
Stacy \cite{Stacy1962} as the class of probability distributions
simultaneously containing both Weibull and gamma distributions. The
family of GG distributions contains practically all the most popular
absolutely continuous distributions concentrated on the non-negative
half-line. In particular, the family of GG distributions contains:

\begin{itemize}

\vspace{-3mm}\item[$\bullet$] the gamma distribution $(\alpha=1)$
and its special cases --

\begin{itemize}

\vspace{-3mm}\item[$\circ$] the exponential distribution
($\alpha=1$, $r =1$),

\vspace{-1mm}\item[$\circ$] the Erlang distribution ($\alpha=1$, $ r
\in\mathbb{N}$),

\vspace{-1mm}\item[$\circ$] the chi-square distribution ($\alpha=1$,
$\lambda=\frac12$);

\end{itemize}

\vspace{-4mm}\item[$\bullet$] the Nakagami distribution
($\alpha=2$);

\vspace{-2mm}\item[$\bullet$] the half-normal (folded normal)
distribution (the distribution of the maximum of a standard Wiener
process on the interval $[0,1]$) ($\alpha=2$, $ r =\frac12$);

\vspace{-2mm}\item[$\bullet$] the Rayleigh distribution ($\alpha=2$,
$ r =1$);

\vspace{-2mm}\item[$\bullet$] the chi-distribution ($\alpha=2$,
$\lambda=1/\sqrt{2}$);

\vspace{-2mm}\item[$\bullet$] the Maxwell distribution (the
distribution of the absolute values of the velocities of moleculas
in a dilute gas) ($\alpha=2$, $ r =\frac32$);

\vspace{-2mm}\item[$\bullet$] the Weibull--Gnedenko distribution
(the extreme value distribution of type III) ($r =1$, $\alpha>0$);

\vspace{-2mm}\item[$\bullet$] the (folded) exponential power
distribution ($\alpha>0$, $r=\frac{1}{\alpha}$);

\vspace{-2mm}\item[$\bullet$] the inverse gamma distribution
($\alpha=-1$) and its special case --

\begin{itemize}

\vspace{-3mm}\item[$\circ$] the L{\'e}vy distribution (the one-sided
stable distribution with the characteristic exponent $\frac12$ --
the distribution of the first hit time of the unit level by the
Brownian motion) ($\alpha=-1$, $r =\frac12$);

\end{itemize}

\vspace{-2mm}\item[$\bullet$] the Fr{\'e}chet distribution (the
extreme value distribution of type II) ($r=1$, $\alpha<0$)

\end{itemize}

\vspace{-2mm} \noindent and other laws. The limit point of the class
of GG distributions is

\begin{itemize}

\vspace{-2mm}\item[$\bullet$] the log-normal distribution ($r
\to\infty$).

\end{itemize}

\vspace{-2mm}

GG distributions are widely applied in many practical problems,
first of all, related to image or signal processing. There are
dozens of papers dealing with the application of GG distributions as
models of regularities observed in practice. Apparently, the
popularity of GG distributions is due to that most of them can serve
as adequate asymptotic approximations, since all the representatives
of the class of GG distributions listed above appear as limit laws
in various limit theorems of probability theory in rather simple
limit schemes. Below we will formulate a general limit theorem (an
analog of the law of large numbers) for random sums of independent
r.v.'s in which the GG distributions are limit laws.

Our interest to these distribution was directly motivated by the
problem of modeling some statistical regularities in precipitation.
In most papers dealing with the statistical analysis of
meteorological data available to the authors, the suggested
analytical models for the observed statistical regularities in
precipitation are rather ideal and far from being adequate. For
example, it is traditionally assumed that the duration of a wet
period (the number of subsequent wet days) follows the geometric
distribution (for example, see~\cite{Zolina2013}) despite that the
goodness-of-fit of this model is very far from being admissible.
Perhaps, this prejudice is based on the conventional interpretation
of the geometric distribution in terms of the Bernoulli trials as
the distribution of the number of subsequent wet days
(``successes'') till the first dry day (``failure''). But the
framework of Bernoulli trials assumes that the trials are
independent whereas a thorough statistical analysis of precipitation
data registered at different points demonstrates that the sequence
of dry and wet days is not only devoid of the independence property,
but is also not Markovian, so that the framework of Bernoulli trials
is absolutely inadequate for analyzing meteorological data.

It turned out that the statistical regularities of the number of
subsequent wet days can be very reliably modeled by the negative
binomial distribution with the shape parameter less than one. For
example, in \cite{Gulev} the data registered in so climatically
different points as Potsdam (Brandenburg, Germany) and Elista
(Kalmykia, Russia) was analyzed and it was demonstrated that the
fluctuations of the numbers of successive wet days with very high
confidence fit the negative binomial distribution with shape
parameter $r\approx 0.85$. In the same paper a schematic attempt was
undertaken to explain this phenomenon by the fact that negative
binomial distributions can be represented as mixed Poisson laws with
mixing gamma-distributions whereas the Poisson distribution is the
best model for the discrete stochastic chaos (see, e. g.,
\cite{Kingman1993, KorolevBeningShorgin2011}) and the mixing
distribution accumulates the stochastic influence of factors that
can be assumed exogenous with respect to the local system under
consideration.

Negative binomial distributions are special cases of generalized
negative binomial (GNB) distributions which are mixed Poisson
distributions with mixing conducted with respect to a GG
distribution. This family of discrete distributions is very wide and
embraces Poisson distributions (as limit points corresponding to a
degenerate mixing distribution), negative binomial (Polya)
distributions including geometric distributions (corresponding to
the gamma mixing distribution, see \cite{GreenwoodYule1920}), Sichel
distributions (corresponding to the inverse gamma mixing
distributions, see \cite{Sichel1971}), Weibull--Poisson
distributions (corresponding to the Weibull mixing distributions,
see \cite{KorolevPoisson}) and many other types supplying
descriptive statistics with many flexible models. More examples of
mixed Poisson laws can be found in \cite{Grandell1997,
SteutelvanHarn2004}.

It is quite natural to expect that, having introduced one more free
parameter into the pure negative binomial model, namely, the power
parameter in the exponent of the original gamma mixing distribution,
instead of the negative binomial model one might obtain a more
flexible GNB model that provides even better fit with the
statistical data of the durations of wet days. The analysis of the
real data shows that this is indeed so. On Fig. 1 there are the
histogram constructed from real data of 3320 wet periods in Potsdam
and the fitted negative binomial distribution (that is, the mixed
Poisson distribution with the GG mixing law with $\alpha=1$) (right)
and the fitted GNB distribution with $\alpha=0.917$ (left). The
$L_1$-distance between the histogram and the fitted GNB distribution
is 0.02637 that is nearly 1.26 times less than the $L_1$-distance
between the histogram and the fitted negative binomial model (equal
to 0.03318).

\renewcommand{\figurename}{{\small\bf Fig.}}

\begin{figure}[h]
\begin{minipage}[h]{0.49\textwidth}
\center{\includegraphics[height=0.8\textwidth, width=\textwidth]{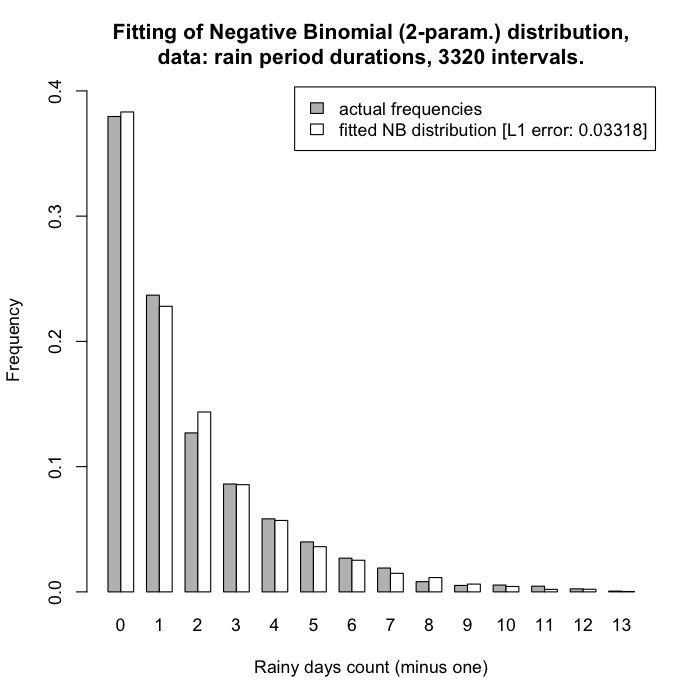} \\
(a)}
\end{minipage}
\hfill
\begin{minipage}[h]{0.49\textwidth}
\center{\includegraphics[height=0.8\textwidth, width=\textwidth]{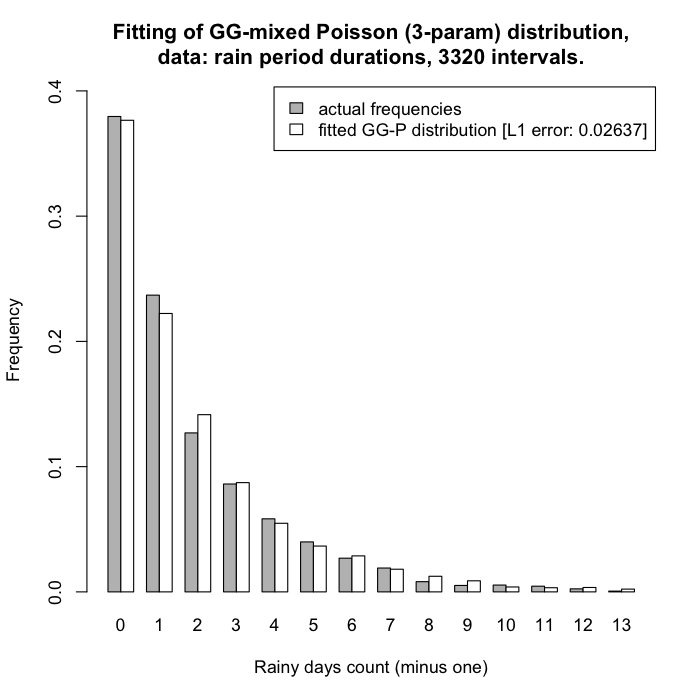} \\
(b)}
\end{minipage}
\label{Data} \caption{{\small The histogram constructed from real
data of 3320 wet periods in Potsdam and the fitted negative binomial
distribution (a) and the fitted GG-mixed Poisson distribution with
$\alpha=0.917$ (b).}}
\end{figure}

The purpose of the present paper is two-fold. First, we try to give
further theoretic explanation of the high adequacy of the negative
binomial model for the description of meteorological data. For this
purpose we use the concept of a mixed geometric law introduced in
\cite{Korolev2016TVP} (also see \cite{KorolevPoisson, Korolev2016}).
Having first proved that any GG distribution with shape parameter
$r$ and exponent power parameter $\gamma$ less than one is mixed
exponential and thus generalizing Gleser's similar theorem on
gamma-distributions \cite{Gleser1989}, we then prove that any mixed
Poisson distribution with the generalized gamma mixing law ({\it GNB
distribution}) with shape parameter and exponent power parameter
less than one is actually mixed geometric. GG distributions with
exponent power less than one are of special interest since they
occupy an intermediate position between distributions with
exponentially decreasing tails (e. g., exponential and
gamma-distributions) and heavy-tailed distributions with
Zipf--Pareto power-type decrease of tails.

The mixed geometric distribution can be interpreted in terms of the
Bernoulli trials as follows. First, as a result of some
<<preliminary>> experiment the value of some random variable taking
values in $[0,1]$ is determined which is then used as the
probability of success in the sequence of Bernoulli trials in which
the original <<unconditional>> mixed Poisson random variable is
nothing else than the <<conditionally>> geometrically distributed
random variable having the sense of the number of trials up to the
first failure. This makes it possible to assume that the sequence of
wet/dry days is not independent, but is conditionally independent
and the random probability of success is determined by some outer
stochastic factors. As such, we can consider the seasonality or the
type of the cause of a rainy period. So, since the GG distribution
is a more general and hence, more flexible model than the ``pure''
gamma-distribution, there arises a hope that the GNB distribution
could provide even better goodness of fit to the statistical
regularities in the duration of wet periods than the ``pure''
negative binomial distribution.

Second, we study analytic and asymptotic properties of probability
models closely related to GNB distributions such as the
distributions of statistics constructed from samples with random
sizes having the GNB distribution. The obtained results can serve as
a theoretical explanation of some mixed models used in the
statistical analysis.

\subsection{Notation and definitions}

In the paper, conventional notation is used. The symbols $\eqd$ and
$\Longrightarrow$ denote the coincidence of distributions and
convergence in distribution, respectively.

In what follows, for brevity and convenience, the results will be
presented in terms of r.v.'s with the corresponding distributions.
It will be assumed that all the r.v.'s are defined on the same
probability space $(\Omega,\,\mathfrak{F},\,{\sf P})$.

A r.v. having the gamma distribution with shape parameter $r>0$ and
scale parameter $\lambda>0$ will be denoted $G_{r,\lambda}$,
$$
{\sf P}(G_{r,\lambda}<x)=\int_{0}^{x}g(z;r,\lambda)dz,\ \
\text{with}\ \
g(x;r,\lambda)=\frac{\lambda^r}{\Gamma(r)}x^{r-1}e^{-\lambda x},\
x\ge0,
$$
where $\Gamma(r)$ is Euler's gamma-function,
$\Gamma(r)=\int_{0}^{\infty}x^{r-1}e^{-x}dx$, $r>0$.

In this notation, obviously, $G_{1,1}$ is a r.v. with the standard
exponential distribution: ${\sf P}(G_{1,1}<x)=\big[1-e^{-x}\big]{\bf
1}(x\ge0)$ (here and in what follows ${\bf 1}(A)$ is the indicator
function of a set $A$).

The gamma distribution is a particular representative of the class
of generalized gamma distributions (GG distributions), which were
first described in \cite{Stacy1962} as a special family of lifetime
distributions containing both gamma and Weibull distributions.

\smallskip

{\sc Definition 1.} A {\it generalized gamma $($GG$)$ distribution}
is the absolutely continuous distribution defined by the density
$$
g^*(x;r,\alpha,\lambda)=\frac{|\alpha|\lambda^r}{\Gamma(r)}x^{\alpha
r-1}e^{-\lambda x^{\alpha}},\ \ \ \ x\ge0,
$$
with $\alpha\in\mathbb{R}$, $\lambda>0$, $r>0$.

\smallskip

The properties of GG distributions are described in \cite{Stacy1962,
KorolevZaks2013}. A r.v. with the density $g^*(x;r,\alpha,\lambda)$
will be denoted $G^*_{r,\alpha,\lambda}$.

For a r.v. with the {\it Weibull distribution}, a particular case of
GG distributions corresponding to the density $g^*(x;1,\alpha,1)$
and the distribution function (d.f.)
$\big[1-e^{-x^{\alpha}}\big]{\bf 1}(x\ge0)$ with $\alpha>0$, we will
use a special notation $W_{\alpha}$. Thus, $G_{1,1}\eqd W_1$. The
density $g^*(x;1,\alpha,1)$ with $\alpha<0$ defines the {\it
Fr{\'e}chet} or {\it inverse Weibull distribution}. It is easy to
see that $W_1^{1/\alpha}\eqd W_{\alpha}$.

In what follows we will be mostly interested in GG distributions
with $\alpha\in(0,1]$.

A r.v. with the standard normal d.f. $\Phi(x)$ will be denoted $X$,
$$
{\sf
P}(X<x)=\Phi(x)=\frac{1}{\sqrt{2\pi}}\int_{-\infty}^{x}e^{-z^2/2}dz,\
\ \ \ x\in\mathbb{R}.
$$

The d.f. and the density of a strictly stable distribution with the
characteristic exponent $\alpha$ and shape parameter $\theta$
defined by the characteristic function (ch.f.)
$$
\mathfrak{f}(t;\alpha,\theta)=\exp\big\{-|t|^{\alpha}\exp\{-{\textstyle\frac12}i\pi\theta\alpha\mathrm{sign}t\}\big\},\
\ \ \ t\in\r,
$$
where $0<\alpha\le2$, $|\theta|\le\min\{1,\frac{2}{\alpha}-1\}$,
will be respectively denoted $F(x;\alpha,\theta)$ and
$f(x;\alpha,\theta)$ (see, e. g., \cite{Zolotarev1983}). A r.v. with
the d.f. $F(x;\alpha,\theta)$ will be denoted $S_{\alpha,\theta}$.
To symmetric strictly stable distributions there correspond the
value $\theta=0$ and the ch.f.
$\mathfrak{f}(t;\alpha,0)=e^{-|t|^{\alpha}}$, $t\in\r$. Hence, it is
easy to see that $S_{2,0}\eqd\sqrt{2}X$.

To one-sided strictly stable distributions concentrated on the
nonnegative halfline there correspond the values $\theta=1$ and
$0<\alpha\le1$. The pairs $\alpha=1$, $\theta=\pm1$ correspond to
the distributions degenerate in $\pm1$, respectively. All the other
strictly stable distributions are absolutely continuous. Stable
densities cannot explicitly be represented via elementary functions
with four exceptions: the normal distribution ($\alpha=2$,
$\theta=0$), the Cauchy distribution ($\alpha=1$, $\theta=0$), the
L{\'e}vy distribution ($\alpha=\frac12$, $\theta=1$) and the
distribution symmetric to the L{\'e}vy law ($\alpha=\frac12$,
$\theta=-1$). Expressions of stable densities in terms of the Fox
functions (generalized Meijer G-functions) can be found in
\cite{Schneider1986, UchaikinZolotarev1999}.

According to the <<multiplication theorem>> (see, e. g., theorem
3.3.1 in \cite{Zolotarev1983}) for any admissible pair of parameters
$(\alpha,\,\theta)$ and any $\alpha'\in(0,1]$ the multiplicative
representation $S_{\alpha\alpha',\theta}\eqd S_{\alpha,\theta}\cdot
S_{\alpha',1}^{1/\alpha}$ holds, in which the factors on the
right-hand side are independent. In particular, for any
$\alpha\in(0,2]$
$$
S_{\alpha,0}\eqd X\sqrt{2S_{\alpha/2,1}},\eqno(1)
$$
that is, any symmetric strictly stable distribution is a normal
scale mixture.

It is well known that if $0<\alpha<2$, then ${\sf
E}|S_{\alpha,\theta}|^{\beta}<\infty$ for any $\beta\in(0,\alpha)$,
but the moments of the r.v. $S_{\alpha,\theta}$ of orders
$\beta\ge\alpha$ do not exist (see, e. g., \cite{Zolotarev1983}).
Despite the absence of explicit expressions for the densities of
stable distributions in terms of elementary functions, it can be
shown \cite{KorolevWeibull2016} that
$$
{\sf E}|S_{\alpha,0}|^{\beta}=\frac{2^{\beta}}{\sqrt{\pi}}\cdot
\frac{\Gamma\big(\frac{\beta+1}{2}\big)\Gamma\big(\frac{\alpha-\beta}{\alpha}\big)}{\Gamma\big(\frac{2-\beta}{\beta}\big)}
$$
for $0<\beta<\alpha<2$ and
$$
{\sf
E}S_{\alpha,1}^{\beta}=\frac{\Gamma\big(1-\frac{\beta}{\alpha}\big)}{\Gamma(1-\beta)}
$$
for $0<\beta<\alpha\le 1$.

In \cite{KotzOstrovskii1996, KorolevZeifman2016a,
KorolevZeifman2016b} it was proved that if $\alpha\in(0,1)$ and the
i.i.d. r.v.'s $S_{\alpha,1}$ and $S'_{\alpha,1}$ have the same
strictly stable distribution, then the density $v_{\alpha}(x)$ of
the r.v. $R_{\alpha}=S_{\alpha,1}/S'_{\alpha,1}$ has the form
$$
v_{\alpha}(x)=\frac{\sin(\pi\alpha)x^{\alpha-1}}{\pi[1+x^{2\alpha}+2x^{\alpha}\cos(\pi\alpha)]},\
\ \ x>0.\eqno(2)
$$

The distribution of a r.v. $Z$ is said to belong to the domain of
normal attraction of the strictly stable law $F(x;\alpha,\theta)$,
$\mathcal{L}(Z)\in DNA\big(F(x;\alpha,\theta)\big)$, if there exists
a finite positive constant $c$ such that
$$
\frac{c}{n^{1/\alpha}}\sum\nolimits_{j=1}^nZ_j\Longrightarrow
S_{\alpha,\theta}\ \ \ (n\to\infty),
$$
where $Z_1,Z_2,...$ are independent copies of the r.v. $Z$. In what
follows we will consider the standard scale case and let $c=1$. In
\cite{Tucker1975} it was shown that if $\mathcal{L}(Z)\in
DNA\big(F(x;\alpha,\theta)\big)$, then ${\sf E}|Z|^{\beta}=\infty$
for any $\beta>\alpha$.

A r.v. $N_{r,p}$ is said to have the {\it negative binomial
distribution} with parameters $r>0$ and $p\in(0,1)$, if
$$
{\sf P}(N_{r,p}=k)=\frac{\Gamma(r+k)}{k!\Gamma(r)}\cdot p^r(1-p)^k,\
\ \ \ k=0,1,2,...
$$

A particular case of the negative binomial distribution
corresponding to the value $r=1$ is the {\it geometric
distribution}. Let $p\in(0,1)$ and let $N_{1,p}$ be the r.v. having
the geometric distribution with parameter $p\,$:
$$
{\sf P}(N_{1,p}=k)=p(1-p)^{k},\ \ \ \ k=0,1,2,...\eqno(3)
$$
This means that for any $m\in\mathbb{N}$
$$
{\sf P}(N_{1,p}\ge
m)=\sum\nolimits_{k=m}^{\infty}p(1-p)^{k}=(1-p)^{m}.
$$

\smallskip

{\sc Definition 2.} Let $Y$ be a r.v. taking values in the interval
$(0,1)$. Moreover, let for all $p\in(0,1)$ the r.v. $Y$ and the
geometrically distributed r.v. $N_{1,p}$ be independent. Let
$V=N_{1,Y}$, that is, $V(\omega)=N_{1,Y(\omega)}(\omega)$ for any
$\omega\in\Omega$. The distribution
$$
{\sf P}(V\ge m)=\int_{0}^{1}(1-y)^{m}d{\sf P}(Y<y), \ \ \
m\in\mathbb{N},
$$
of the r.v. $V$ will be called {\it $Y$-mixed geometric}
\cite{Korolev2016TVP}.

\smallskip

Let $\alpha\in(0,1]$. The distribution with the Laplace--Stieltjes
transform (L.--S.t.)
$$
\mathfrak{m}(s;\alpha)=\frac{1}{1+s^{\alpha}},\ \ \ s\ge0.\eqno(4)
$$
is conventionally called the {\it Mittag-Leffler distribution}. The
origin of this term is due to the form of the probability density
$$
f^{M}(x;\alpha)=\frac{1}{x^{1-\alpha}}\sum\nolimits_{n=0}^{\infty}\frac{(-1)^nx^{\alpha
n}}{\Gamma(\alpha n+1)}=-\frac{d}{dx}E_{\alpha}(-x^{\alpha}),\ \ \
x\ge0,\eqno(5)
$$
corresponding to L.--S.t. (4), where $E_{\alpha}(z)$ is the
Mittag-Leffler function with index $\alpha$ that is defined as the
power series
$$
E_{\alpha}(z)=\sum \nolimits_{n=0}^{\infty}\frac{z^n}{\Gamma(\alpha
n+1)},\ \ \ \alpha>0,\ z\in\mathbb{Z}.
$$
The distribution function corresponding to density (5) will be
denoted $F^{M}(x;\alpha)$. The r.v. with the d.f. $F^{M}(x;\alpha)$
will be denoted $M_{\alpha}$. In \cite{KorolevZeifman2016b} the
integral representation
$$
f^{M}(x;\alpha)=\frac{\sin(\pi\alpha)}{\pi}\int_{0}^{\infty}\frac{
z^{\alpha}e^{-zx}dz}{1+z^{2\alpha}+2z^{\alpha}\cos(\pi\alpha)},\ \ \
x>0,
$$
for the Mittag-Leffler density was proved.

With $\alpha=1$, the Mittag-Leffler distribution turns into the
standard exponential distribution, that is, $M_1\eqd W_1$. But with
$\alpha<1$ the Mittag-Leffler distribution density has a heavy
power-type tail: from the well-known asymptotic properties of the
Mittag-Leffler function it can be deduced that if $0<\alpha<1$, then
$$
f^{M}(x;\alpha)\sim \frac{\sin(\alpha\pi)\Gamma(\alpha+1)}{\pi
x^{\alpha+1}}
$$
as $x\to\infty$, see, e. g., \cite{Kilbas2014}.

It is well-known that the Mittag-Leffler distribution is stable with
respect to geometric summation (or {\it geometrically stable}). This
means that if $X_1,X_2,...$ are independent random variables and
$N_{1,p}$ is the random variable independent of $X_1,X_2,...$ and
having the geometric distribution (3), then for each $p\in(0,1)$
there exists a constant $a_p>0$ such that
$a_p\big(X_1+...+X_{N_{1,p}}\big)\Longrightarrow M_{\alpha}$ as
$p\to 0$, see, e. g., \cite{Bunge1996} or \cite{KlebanovRachev1996}.
Moreover, as long ago as in 1965 it was shown by I.~Kovalenko
\cite{Kovalenko1965} that the distributions with L.--S.t:s (4) are
the only possible limit laws for the distributions of appropriately
normalized geometric sums of the form
$a_p\big(X_1+...+X_{N_{1,p}}\big)$ as $p\to0$, where $X_1,X_2,...$
are independent identically distributed {\it nonnegative} random
variables and $N_{1,p}$ is the random variable with geometric
distribution (3) independent of the sequence $X_1,X_2,...$ for each
$p\in(0,1)$. The proofs of this result were reproduced in
\cite{GnedenkoKovalenko1968, GnedenkoKovalenko1989} and
\cite{GnedenkoKorolev1996}. In these books the class of
distributions with Laplace transforms (4) was not identified as the
class of Mittag-Leffler distributions but was called {\it class}
$\mathcal{K}$ after I.~Kovalenko.

Twenty five years later this limit property of the Mittag-Leffler
distributions was re-discovered by A.~Pillai in \cite{Pillai1989,
Pillai1990} who proposed the term {\it Mittag-Leffler distribution}
for the distribution with Laplace transform (1). Perhaps, since the
works \cite{Kovalenko1965, GnedenkoKovalenko1968,
GnedenkoKovalenko1989} were not easily available to probabilists,
the term {\it class $\mathcal{K}$ distribution} did not take roots
in the literature whereas the term {\it Mittag-Leffler distribution}
became conventional.

In \cite{KotzOstrovskii1996, KorolevZeifman2016a,
KorolevZeifman2016b} it was shown that the Mittag-Leffler
distribution is mixed exponential. Namely,
$$
M_{\alpha}\eqd W_1\cdot R_{\alpha},\eqno(6)
$$
where $R_{\alpha}$ is the r.v. with density (2) independent of
$W_1$.

The Mittag-Leffler distributions are of serious theoretical interest
in problems related to thinned (or rarefied) homogeneous flows of
events such as renewal processes or anomalous diffusion or
relaxation phenomena, see \cite{WeronKotulski1996,
GorenfloMainardi2006} and the references therein.

In 1953 Yu. V. Linnik \cite{Linnik1953} introduced the class of
symmetric probability distributions defined by the characteristic
functions
$$
\mathfrak{f}^L(t;\alpha)=\frac{1}{1+|t|^{\alpha}},\ \ \
t\in\mathbb{R},\eqno(7)
$$
where $\alpha\in(0,2]$. Later the distributions of this class were
called {\it Linnik distributions} \cite{Kotz2001} or {\it
$\alpha$-Laplace distributions} \cite{Pillai1985}. In this paper we
will keep to the first term that has become conventional. With
$\alpha=2$, the Linnik distribution turns into the Laplace
distribution corresponding to the density
$$
f^{\Lambda}(x)=\textstyle{\frac12}e^{-|x|},\ \ \
x\in\mathbb{R}.\eqno(8)
$$
A random variable with Laplace density (8) and its distribution
function will be denoted $\Lambda$ and $F^{\Lambda}(x)$,
respectively.

A random variable with the Linnik distribution with parameter
$\alpha$ will be denoted $L_{\alpha}$. Its distribution function and
density will be denoted $F^{L}(x;\alpha)$ and $f^{L}(x;\alpha)$,
respectively. Moreover, from (7) and (8) it follows that
$F^{L}(x;2)\equiv F^{\Lambda}(x)$, $x\in\mathbb{R}$.

The Linnik distributions possess many interesting properties. First
of all, as well as the Mittag-Leffler distributions, the Linnik
distributions are geometrically stable, that is, if $X_1,X_2,...$
are i.i.d. r.v.'s with $\mathcal{L}(X_1)\in
DNA\big(F(x;\alpha,0)\big)$, then under an appropriate choice of
positive constants $a_p$ the distributions of the normalized
geometric random sums $a_p(X_1+...+X_{N_{1,p}})$ converge to the
Linnik distribution with parameter $\alpha$.

The Linnik distributions are unimodal \cite{Laha1961} and infinitely
divisible \cite{Devroye1990}, have an infinite peak of the density
for $\alpha\le1$ \cite{Devroye1990}, etc. In
\cite{KotzOstrovskiiHayfavi1995a, KotzOstrovskiiHayfavi1995b} a
detailed investigation of analytic and asymptotic properties of the
density of the Linnik distribution was carried out. In the papers
\cite{KotzOstrovskii1996, KorolevZeifman2016b} it was demonstrated
that the Linnik distributions are intimately related with the
normal, Laplace, stable and Mittag-Leffler distributions:
$$
L_{\alpha}\eqd X\sqrt{2M_{\alpha/2}}\eqd
\Lambda\sqrt{R_{\alpha/2}},\eqno(9)
$$
where all the factors are independent and the r.v. $R_{\alpha/2}$ is
the ratio of two i.i.d. strictly stable r.v.'s $S_{\alpha/2,1}$ and
$S'_{\alpha/2,1}$ and has the density (2) with $\alpha$ replaced by
$\alpha/2$.

\section{Generalized negative binomial and related distributions}

\subsection{Continuous case. Generalization of Gleser's theorem to
GG distributions}

In the paper \cite{Gleser1989} it was shown that any gamma
distribution with shape parameter no greater than one is mixed
exponential. For convenience, we formulate this result as the
following lemma.

\smallskip

{\sc Lemma 1} \cite{Gleser1989}. {\it The density of a gamma
distribution $g(x;r,\mu)$ with $0<r<1$ can be represented as
$$
g(x;r,\mu)=\int_{0}^{\infty}ze^{-zx}p(z;r,\mu)dz,
$$
where
$$
p(z;r,\mu)=\frac{\mu^r}{\Gamma(1-r)\Gamma(r)}\cdot\frac{\mathbf{1}(z\ge\mu)}{(z-\mu)^rz}.
$$
Moreover, a gamma distribution with shape parameter $r>1$ cannot be
represented as a mixed exponential distribution.}

\smallskip

{\sc Lemma 2} \cite{Korolev2017}. {\it For $r\in(0,1)$ let
$G_{r,\,1}$ and $G_{1-r,\,1}$ be independent gamma-distributed
r.v.'s. Let $\mu>0$. Then the density $p(z;r,\mu)$ in lemma $1$
corresponds to the r.v.
$$
Z_{r,\mu}=\frac{\mu(G_{r,\,1}+G_{1-r,\,1})}{G_{r,\,1}}\eqd\mu
Z_{r,1}\eqd\mu\big(1+\textstyle{\frac{1-r}{r}}Q_{1-r,r}\big),
$$
where $Q_{1-r,r}$ is the r.v. with the Snedecor--Fisher distribution
defined by the probability density
$$
q(x;1-r,r)=\frac{(1-r)^{1-r}r^r}{\Gamma(1-r)\Gamma(r)}
\cdot\frac{1}{x^{r}[r+(1-r)x]},\ \ \ x\ge0.
$$
}

\smallskip

{\sc Remark 1}. It is easily seen that the sum $G_{r,1}+G_{1-r,1}$
has the standard exponential distribution: $G_{r,1}+G_{1-r,1}\eqd
W_1$. However, the numerator and denominator of the expression
defining the r.v. $Z_{r,\mu}$ are not independent.

\smallskip

{\sc Lemma 3} \cite{ShanbhagSreehari1977, KorolevWeibull2016}. {\it
Let $\alpha\in(0,1]$. Then $W_1^{1/\alpha}\eqd W_{\alpha}\eqd
W_1\cdot S_{\alpha,1}^{-1}$ with the r.v.'s on the right-hand side
being independent.}

\smallskip

{\sc Lemma 4}. {\it A d. f. $F(x)$ with $F(0)=0$ corresponds to a
mixed exponential distribution if and only if the function $1-F(x)$
is completely monotone$:$ $F\in C_{\infty}$ and
$(-1)^{n+1}F^{(n)}(x)\ge 0$ for all $x>0$}.

\smallskip

This statement immediately follows from the Bernstein theorem
\cite{Bernstein1928}.

\smallskip

{\sc Theorem 1}. $(i)$ {\it Let $\alpha\in(0,1]$, $r\in(0,1)$,
$\mu>0$. Then the GG distribution with parameters $r$, $\alpha$,
$\mu$ is a mixed exponential distribution$:$ $G^*_{r,\alpha,\mu}\eqd
W_1\cdot\big(S_{\alpha,1}Z_{r,\mu}^{1/\alpha}\big)^{-1}$ with the
r.v.'s on the right-hand side being independent.

\noindent $(ii)$ Let $\alpha>0$, $r>0$, $\mu>0$. A GG distribution
with $\alpha r>1$ cannot be represented as mixed exponential.}

\smallskip

{\sc Proof.} $(i)$. First, note that ${\sf
P}(G_{r,\mu}^{1/\alpha}>x)={\sf P}(G_{r,\mu}>x^{\alpha})$. Hence,
according to lemma 1 for $x\ge0$ we have
$$
{\sf P}(G_{r,\mu}^{1/\alpha}>x)={\sf
P}(W_1>Z_{r,\mu}x^{\alpha})=\int_{0}^{\infty}\!\!\!e^{-zx^{\alpha}}p(z;r,\mu)dz=
\int_{0}^{\infty}\!\!\!{\sf
P}(W_{\alpha}>xz^{1/\alpha})p(z;r,\mu)dz,
$$
that is, $G_{r,\mu}^{1/\alpha}\eqd W_{\alpha}\cdot
Z_{r,\mu}^{-1/\alpha}$. Now apply lemma 3 and obtain
$$
G_{r,\mu}^{1/\alpha}\eqd W_1\cdot
\big(S_{\alpha,1}Z_{r,\mu}^{1/\alpha})^{-1}.\eqno(10)
$$
Second, it is easy to see that
$$
G_{r,\mu}^{1/\alpha}\eqd G^*_{r,\alpha,\mu}\eqno(11)
$$
for any $r>0$, $\mu>0$ and $\alpha>0$. Now the desired assertion
follows from (10) and (11).

To prove assertion $(ii)$, assume that $\alpha r>1$ and the r.v.
$G^*_{r,\alpha,\mu}$ has a mixed exponential distribution. By lemma
4 this means that the function $\psi(s)={\sf
P}(G^*_{r,\alpha,\mu}>s)$, $s\ge0$, is completely monotone. But
$\psi'(s)=g^*(s;r,\alpha,\mu)\ge0$ for all $s\ge0$, whereas
$$
\psi''(s)=(g^*)'(s;r,\alpha,\mu)=\frac{\alpha\mu^r}{\Gamma(r)}s^{\alpha
r-2}e^{-\mu s^{\alpha}}\big((\alpha r-1)-\mu\alpha
s^{\alpha}\big)\le0,
$$
only if $(\alpha r-1)-\mu\alpha s^{\alpha}\le0$, that is, $s\ge
s_0\equiv \big[(\alpha r-1)/\mu\alpha\big]^{1/\alpha}>0$, and
$\psi''(s)\ge0$ for $s\in(0,s_0)\neq\varnothing$ contradicting the
complete monotonicity of $\psi(s)$ and thus proving the second
assertion. The theorem is proved.

\smallskip

{\sc Remark 2.} Lemma 3 states that the Weibull distribution with
parameter $\alpha>0$ is mixed exponential, if $\alpha\le1$. Now from
Theorem 1 it follows that this statement can be reinforced: the
Weibull distribution is mixed exponential, {\it if and only if}
$\alpha\le1$.

\subsection{Discrete case. An analog of Gleser's theorem for
generalized negative binomial distributions}

{\sc Definition 3}. For $r>0$, $\alpha\in\r$ and $\mu>0$ let
$N_{r,\alpha,\mu}$ be a r.v. with the {\it generalized negative
binomial $($GNB$)$ distribution}:
$$
{\sf
P}(N_{r,\alpha,\mu}=k)=\frac{1}{k!}\int_{0}^{\infty}e^{-z}z^kg^*(z;r,\alpha,\mu)dz,\
\ \ \ k=0,1,2...
$$

\smallskip

%

{\sc Theorem 2}. {\it If $r\in(0,1]$, $\alpha\in(0,1]$ and $\mu>0$,
then a GNB distribution is a $Y_{r,\alpha,\mu}$-mixed geometric
distribution$:$
$$
{\sf P}(N_{r,\alpha,\mu}=k)=\int_{0}^{1}y(1-y)^kd{\sf
P}(Y_{r,\alpha,\mu}<y),\ \ \ \ k=0,1,2...,
$$
where
$$
Y_{r,\alpha,\mu}\eqd
\frac{S_{\alpha,1}Z_{r,\mu}^{1/\alpha}}{1+S_{\alpha,1}Z_{r,\mu}^{1/\alpha}}\eqd
\frac{\mu^{1/\alpha}S_{\alpha,1}(G_{r,1}+G_{1-r,1})^{1/\alpha}}{G_{r,1}^{1/\alpha}+
\mu^{1/\alpha}S_{\alpha,1}(G_{r,1}+G_{1-r,1})^{1/\alpha}},\eqno(12)
$$
where the r.v.'s $S_{\alpha,1}$ and $Z_{\mu,r}$ or $S_{\alpha,1}$,
$G_{r,1}$ and $G_{1-r,1}$ are independent.}

\smallskip

{\sc Proof}. Using theorem 1 we have
$$
{\sf
P}(N_{r,\alpha,\mu}=k)=-\frac{1}{k!}\int_{0}^{\infty}e^{-z}z^kd{\sf
P}(G^*_{r,\alpha,\mu}>z)=-\frac{1}{k!}\int_{0}^{\infty}e^{-z}z^kd{\sf
P}(W_1>S_{\alpha,1}Z_{r,\mu}^{1/\alpha}z)=
$$
$$
=\frac{1}{k!}\int_{0}^{\infty}x\bigg(\int_{0}^{\infty}e^{-z(1+x)}z^kdz\bigg)d{\sf
P}(S_{\alpha,1}Z_{r,\mu}^{1/\alpha}<x)=
$$
$$
=\frac{\Gamma(k+1)}{k!}\int_{0}^{\infty}\frac{x}{(1+x)^{k+1}}d{\sf
P}(S_{\alpha,1}Z_{r,\mu}^{1/\alpha}<x)=
\int_{0}^{\infty}\frac{x}{1+x}\Big(1-\frac{x}{1+x}\Big)^kd{\sf
P}(S_{\alpha,1}Z_{r,\mu}^{1/\alpha}<x).
$$
Changing the variables $\frac{x}{1+x}\longmapsto y$, we finally
obtain
$$
{\sf P}(N_{r,\alpha,\mu}=k)=\int_{0}^{1}y(1-y)^kd{\sf
P}\Big(S_{\alpha,1}Z_{r,\mu}^{1/\alpha}<\frac{y}{1-y}\Big)=\int_{0}^{1}y(1-y)^kd{\sf
P}\bigg(\frac{S_{\alpha,1}Z_{r,\mu}^{1/\alpha}}{1+S_{\alpha,1}Z_{r,\mu}^{1/\alpha}}<y\bigg).\eqno(13)
$$
Moreover, (13) and lemma 2 yield representation (12). The theorem is
proved.

\smallskip

{\sc Remark 3}. Using lemma 1 it is easy to verify that the density
$q(y;r,\alpha,\mu)$ of the r.v. $Y_{r,\alpha,\mu}$ admits the
following integral representation via the strictly stable density
$f(x;\alpha,1)$:
$$
q(y;r,\alpha,\mu)=\frac{\mu^r}{\Gamma(1-r)\Gamma(r)}\cdot\frac{1}{(1-y)^2}\int_{\mu}^{\infty}
\frac{f\big(y(1-y)^{-1}z^{-1/\alpha};\alpha,1\big)dz}{(z-\mu)^rz^{1+2/\alpha}},\
\ \ 0\le y\le1.
$$

\smallskip

From (12) we easily obtain the following asymptotic assertion.

\smallskip

{\sc Corollary 1}. {\it As $\mu\to0$, the r.v. $Y_{r,\alpha,\mu}$ is
the quantity of order $\mu^{1/\alpha}$ in the sense that
$$
\mu^{-1/\alpha}Y_{r,\alpha,\mu}\Longrightarrow
S_{\alpha,1}Z_{r,1}^{1/\alpha}\eqd
S_{\alpha,1}\cdot\bigg(\frac{G_{r,1}+G_{1-r,1}}{G_{r,1}}\bigg)^{1/\alpha}\eqd
S_{\alpha,1}\cdot\big(1+{\textstyle\frac{1-r}{r}}Q_{1-r,r}\big)^{1/\alpha},\eqno(14)
$$
where the r.v.'s $S_{\alpha,1}$ and $Z_{\mu,r}$ or $S_{\alpha,1}$,
$G_{r,1}$ and $G_{1-r,1}$ are independent and the r.v. $Q_{1-r,r}$
has the Snedecor--Fisher distribution introduced in Lemma $2$.}

\smallskip

{\sc Theorem 3}. {\it Let $r>0$, $\alpha\in\mathbb{R}$, $\mu>0$. We
have
$$
\mu^{1/\alpha}N_{r,\alpha,\mu}\Longrightarrow
G^*_{r,\alpha,1}\eqno(15)
$$
as $\mu\to0$. If, moreover, $r\in(0,1]$ and $\alpha\in(0,1]$, then
the limit law can be represented as
$$
G^*_{r,\alpha,1}\eqd
\frac{W_1}{S_{\alpha,1}Z_{r,1}^{1/\alpha}}\eqd\frac{W_{\alpha}}{Z_{r,1}^{1/\alpha}}\eqd
\bigg(\frac{W_1G_{r,1}}{G_{r,1}+G_{1-r,1}}\bigg)^{1/\alpha}\eqd
W_{\alpha}\cdot\big(1+{\textstyle\frac{1-r}{r}}Q_{1-r,r}\big)^{-1/\alpha},\eqno(16)
$$
where the r.v.'s $W_1$, $S_{\alpha,1}$ and $Z_{r,1}$ are independent
as well as the r.v.'s $W_{\alpha}$ and $Z_{r,1}$, or the r.v.'s
$W_{\alpha}$, $G_{r,1}$ and $G_{1-r,1}$, and the r.v. $Q_{1-r,r}$
has the Snedecor--Fisher distribution introduced in Lemma $2$.}

\smallskip

{\sc Proof}. First, $G^*_{r,\alpha,\mu}\eqd G_{r,\mu}^{1/\alpha}\eqd
\mu^{-1/\alpha}G_{r,1}^{1/\alpha}\eqd\mu^{-1/\alpha}G^*_{r,\alpha,1}$.
Therefore, $\mu^{1/\alpha}G^*_{r,\alpha,\mu}\eqd G^*_{r,\alpha,1}$,
$\mu>0$. Let $P(t)$, $t\ge0$, be a standard Poisson process
(homogeneous Poisson process with unit intensity). Then we obviously
have $N_{r,\alpha,\mu}\eqd P(G^*_{r,\alpha,\mu})$ where the r.v.
$G^*_{r,\alpha,\mu}$ is independent of the process $P(t)$. Hence, by
virtue of lemma 2 of \cite{Korolev1998}, we have (15). Second,
representation (16) follows from theorem 1, corollary 1 and lemma 3.

\smallskip

{\sc Remark 4.} In the case $r\in(0,1]$ and $\alpha\in(0,1]$ the
convergence
$$
\mu^{1/\alpha}N_{r,\alpha,\mu}\Longrightarrow
\frac{W_1}{S_{\alpha,1}Z_{r,1}^{1/\alpha}}\eqd\frac{W_{\alpha}}{Z_{r,1}^{1/\alpha}}\eqd
\bigg(\frac{W_1G_{r,1}}{G_{r,1}+G_{1-r,1}}\bigg)^{1/\alpha}\eqd
W_{\alpha}\cdot\big(1+{\textstyle\frac{1-r}{r}Q_{1-r,r}}\big)^{-1/\alpha}
$$
as $\mu\to 0$ can be obtained as a simple corollary of theorem 2 and
theorem 1 of \cite{Korolev2016TVP} establishing the conditions for
the convergence of mixed geometric distributions.

\smallskip

{\sc Remark 5.} Let $r>0$, $p\in(0,1)$. It is easy to see that
$N_{r,p}\eqd N_{r,1,p/(1-p)}$. Therefore, to obtain the
corresponding results for the negative binomial distribution one
should just set $\alpha=1$, $\mu=p/(1-p)$ and let $p\to0$ in the
above statements under the convention that $S_{1,1}\eqd 1$.

\subsection{Mixed binomial distributions and an analog
of the Poisson theorem}

Consider one more problem related to the scheme of Bernoulli trials
with a random probability $Y_{r,\alpha,\mu}$ of success under the
assumption that this probability is infinitesimal. Within the
framework of this scheme, first, the value of the r.v.
$Y_{r,\alpha,\mu}(\omega)\in(0,1)$ is determined as a result of the
<<preliminary>> experiment. Then this value is assigned to the
probability of success in the sequence of $m\in\mathbb{N}$ Bernoulli
trials. Then the r.v. $K=K(\omega)$ is determined as the number of
successes in $m$ Bernoulli trials with the probability of success
equal to $Y_{r,\alpha,\mu}(\omega)$. To formalize the
infinitesimality of the random probability of success
$Y_{r,\alpha,\mu}$, supply the parameters $\mu$ and $m$ (for the
sake of generality), and, correspondingly, the r.v. $K$ with an
<<infinitely large>> index $n$, which makes it possible to trace the
convergence of the sequence of the r.v.'s
$\{Y_{r,\alpha,\mu_n}\}_{n\ge1}$ to zero as $n\to\infty$. In turn,
the infinitesimality of $Y_{r,\alpha,\mu_n}$ means that successes
are rare events within the scheme of Bernoulli trials with a random
probability of success under consideration.

In the papers \cite{KorolevPoisson, Korolev2016} a <<random>>
version of the classical Poisson theorem (the so-called <<law of
small numbers>>) for {\it mixed binomial distributions} with a
random probability of success and infinitely increasing integer
parameter $m_n$ (<<the number of trials>>) was considered. In the
preceding papers dealing with <<random>> versions of the Poisson
theorem (see, e. g., \cite{KorolevBeningShorgin2011}), on the
contrary, the number of trials was random whereas the probability of
success remained non-random.

Let $Y$ be a r.v. such that ${\sf P}(0<Y<1)=1$, $m\in\mathbb{N}$,
$k=1,2,...$ We will say that the r.v. $K$ has the {\it $Y$-mixed
binomial distribution} with the parameter $m$, if
$$ {\sf
P}(K=j)=\Big(\begin{array}{c}\!\!m\!\!\\ \!\!j\!\!\end{array}\Big)
\int_{0}^{1}z^j(1-z)^{m-j}d{\sf P}(Y<z),\ \ \ j=0,1,...,m.\eqno(17)
$$

Now for each $n\in\mathbb{N}$, let $Y_n$ be a r.v. such that ${\sf
P}(0<Y_n<1)=1$, $m_n\in\mathbb{N}$ and $K_n$ be a r.v. with the
$Y_n$-mixed binomial distribution. For $x\in\mathbb{R}$ denote
$B_n(x)={\sf P}(K_n<x)$. Let $Z$ be a positive r.v. The mixed
Poisson d.f. with the structural r.v. $Z$ (in the terminology of
\cite{Grandell1997}) will be denoted $\Pi^{(Z)}(x)$:
$$
\Pi^{(Z)}(x+0)=\sum\nolimits_{j=0}^{[x]}\frac{1}{j!}\int_{0}^{\infty}e^{-z}z^jd{\sf
P}(Z<z),\ \ \ \ \ \ x\in\mathbb{R}.
$$

\smallskip

{\sc Theorem 4} \cite{Korolev2016}. {\it Let $\{m_n\}_{n\ge1}$ be an
infinitely increasing sequence of natural numbers. For any
$n\in\mathbb{N}$, let $K_n$ be a r.v. with the $Y_n$-mixed binomial
distribution $(12)$ with the parameter $m_n$ and d.f. $B_n(x)$.
Assume that in $(12)$ the r.v.'s $Y_n$ are infinitesimal in the
sense that there exists a r.v. $Z$ such that ${\sf P}(0<
Z<\infty)=1$ and
$$
m_nY_n\Longrightarrow Z \eqno(18)
$$
as $n\to\infty$. Then}
$$
B_n(x)\Longrightarrow \Pi^{(Z)}(x) \ \ \ \ \ (n\to\infty).
$$

\smallskip

Let $r\in(0,1)$, $\alpha\in(0,1)$ and $\{m_n\}_{n\ge1}$ be an
infinitely increasing sequence of natural numbers. For each
$n\in\mathbb{N}$ consider the r.v. $Y_{r,\alpha,1/m_n^{\alpha}}$.
Then from corollary 1 it follows that
$$
m_nY_{r,\alpha,1/m_n^{\alpha}}\Longrightarrow
S_{\alpha,1}Z_{r,1}^{1/\alpha}\eqd
S_{\alpha,1}\cdot\bigg(\frac{G_{r,1}+G_{1-r,1}}{G_{r,1}}\bigg)^{1/\alpha}\eqd
S_{\alpha,1}\big(1+{\textstyle\frac{1-r}{r}}Q_{1-r,r}\big)^{1/\alpha}
\eqno(19)
$$
as $n\to\infty$, where the r.v.'s $S_{\alpha,1}$ and $Z_{\mu,r}$ or
$S_{\alpha,1}$, $G_{r,1}$ and $G_{1-r,1}$ are independent. This
means that condition (18) holds with
$Z=S_{\alpha,1}Z_{r,1}^{1/\alpha}$. As this is so, from theorem 4 we
immediately obtain the following random version of the Poisson
theorem.

\smallskip

{\sc Corollary 3}. {\it Let $r\in(0,1)$, $\alpha\in(0,1)$ and
$\{m_n\}_{n\ge1}$ be an infinitely increasing sequence of natural
numbers. Let for each $n\in\mathbb{N}$ the r.v.
$K_{r,\alpha,1/m_n^{\alpha}}$ have the
$Y_{r,\alpha,1/m_n^{\alpha}}$-mixed binomial distribution with
parameter $m_n$. Then}
$$
{\sf P}(K_{r,\alpha,1/m_n^{\alpha}}<x)\Longrightarrow
\Pi^{(S_{\alpha,1}Z_{r,1}^{1/\alpha})}(x)\ \ \ \ (n\to\infty).
$$

\smallskip

If $P_{r,\alpha}$ is the r.v. with the d.f.
$\Pi^{(S_{\alpha,1}Z_{r,1}^{1/\alpha})}(x)$, then for any
$k=0,1,...$ we have
$$
{\sf P}(P_{r,\alpha}=k)=\frac{1}{k!}\int_{0}^{\infty}z^ke^{-z}d{\sf
P}(S_{\alpha,1}Z_{r,1}^{1/\alpha}<z)= \frac{1}{k!}{\sf
E}\big(S_{\alpha,1}^kZ_{r,1}^{k/\alpha}\exp\{-S_{\alpha,1}Z_{r,1}^{1/\alpha}\}\big).
$$

\section{Limit theorems for sums of independent random variables in
which the number of summands has the GNB distributions}

Some results presented below will substantially rely on the
following auxiliary statement. Consider a sequence of r.v.'s
$Q_1,Q_2,...$ Let $N_1,N_2,...$ be natural-valued r.v.'s such that
for every $n\in\mathbb{N}$ the r.v. $N_n$ is independent of the
sequence $Q_1,Q_2,...$ In the following statement the convergence is
meant as $n\to\infty$.

\smallskip

{\sc Lemma 5.} {\it Assume that there exist an infinitely increasing
$($convergent to zero$)$ sequence of positive numbers
$\{b_n\}_{n\ge1}$ and a r.v. $Q$ such that
$$
b_n^{-1}Q_n\Longrightarrow Q.
$$
If there exist an infinitely increasing $($convergent to zero$)$
sequence of positive numbers $\{d_n\}_{n\ge1}$ and a r.v. $V$ such
that
$$
d_n^{-1}b_{N_n}\Longrightarrow V,\eqno(20)
$$
then
$$
d_n^{-1}Q_{N_n}\Longrightarrow Q\cdot V,\eqno(21)
$$
where the r.v.'s on the right-hand side of $(3)$ are independent.
If, in addition, $N_n\longrightarrow\infty$ in probability and the
family of scale mixtures of the d.f. of the r.v. $Q$ is
identifiable, then condition $(20)$ is not only sufficient for
$(21)$, but is necessary as well.}

\smallskip

{\sc Proof} see in \cite{Korolev1994} (the case $b_n,d_n\to\infty$),
\cite{Korolev1995} (the case $b_n,d_n\to 0$) or
\cite{BeningKorolev2002}, theorem 3.5.5.

\subsection{An analog of the law of large numbers for nonnegative summands.
Generalized R{\'e}nyi theorem}

Consider i.i.d. nonnegative r.v.'s $X_1,X_2,...$ For
$k\in\mathbb{N}$ denote $\Sigma_k=X_1+...+X_k$. Let for each
$n\in\mathbb{N}$ the r.v. $N_{r,\alpha,\mu}$ have the GNB
distribution (see definition 3).

In what follows we will be interested in the asymptotic behavior of
GNB random sums $\Sigma_{N_{r,\alpha,\mu}}$ as $\mu\to 0$.

We begin with the situation where ${\sf E}X_1\equiv c\in(0,\infty)$.
In accordance with Kolmogorov's strong law of large numbers this
condition is in some sense equivalent to that
$$
\frac{1}{n}\sum\nolimits_{j=1}^nX_j\longrightarrow c \eqno(22)  
$$
almost surely as $n\to\infty$. Our nearest aim is to obtain an
analog of the law of large numbers for GNB random sums
$\Sigma_{N_{r,\alpha,\mu}}$ as $\mu\to 0$. Let $\mu=1/n^{\alpha}$.
Consider the r.v. $N_{r,\alpha,1/n^{\alpha}}$.

\smallskip

{\sc Theorem 5.} {\it Assume that i.i.d. nonnegative r.v.'s
$X_1,X_2,...$ satisfy condition $(22)$. Let for each
$n\in\mathbb{N}$ the r.v. $N_{r,\alpha,1/n^{\alpha}}$ have the GNB
distribution with parameters $r>0$, $\alpha\in\mathbb{R}$,
$\mu=1/n^{\alpha}$ and be independent of the sequence $X_1,X_2,...$
Then
$$
\lim_{n\to\infty}\sup_{x\ge0}\bigg|{\sf
P}(\Sigma_{N_{r,\alpha,1/n^{\alpha}}}<nx)-\int_{0}^{x}g^*(z;r,\alpha,1/c^{\alpha})dz\bigg|=0.\eqno(23)
$$
If, moreover, $r\in(0,1)$ and $\alpha\in(0,1)$, then the limit GG
distribution is mixed exponential$:$}
$$
G^*_{r,\alpha,1/c^{\alpha}}\eqd\frac{cW_1}{S_{\alpha,1}Z_{r,1}^{1/\alpha}}.\eqno(24)
$$

\smallskip

{\sc Proof.} From condition (22) it follows that in lemma 5 we can
take $b_n\equiv n$, $Q_n\equiv \Sigma_n$, $Q\equiv c$. From theorem
3 it follows that
$$
n^{-1}N_{r,\alpha,1/n^{\alpha}}\Longrightarrow G^*_{r,\alpha,1},
\eqno(25)
$$
so that we can take $d_n\equiv n$, $V\eqd G^*_{r,\alpha,1}$. By
virtue of lemma 5 conditions (22) and (25) imply
$$
n^{-1}\Sigma_{N_{r,\alpha,1/n^{\alpha}}}\Longrightarrow c
G^*_{r,\alpha,1}\eqd G^*_{r,\alpha,1/c^{\alpha}}\eqno(26)
$$
as $n\to\infty$. The observation that the limit GG distribution is
absolutely continuous leads to the conclusion that the convergence
in distribution (26) is uniform over $x\ge0$. Relation (23) is thus
proved.

In the case $r\in(0,1)$ and $\alpha\in(0,1)$ representation (24)
follows from theorem 1. The theorem is proved.

\smallskip

If $r\in(0,1)$ and $\alpha\in(0,1)$, then theorem 5 can be regarded
as a generalization of the R{\'e}nyi theorem on the asymptotic
behavior of rarefied (thinned) stationary point processes (see, e.
g., \cite{Kalashnikov1997}). The R{\'e}nyi theorem can be regarded
as a law of large numbers for geometric sums. The classical
R{\'e}nyi theorem establishes that a stationary point process
converges to the Poisson process under the ordinary rarefaction when
each point is deleted with probability $1-p$ and left as it is with
probability $p$ accompanied by an appropriate change of scale to
provide the non-degenerateness of the limit process. As is known,
the Poisson process is a renewal process with exponentially
distributed spacings. Theorem 5 generalizes the R{\'e}nyi theorem to
mixed geometric sums. In terms of rarefaction or thinning this can
be interpreted as that the rarefaction becoming <<doubly
stochastic>> so that, prior to rarefaction, the probability $p$ is
chosen at random as $p=Y_{r,\alpha,1/n^{\alpha}}$. This results in
that the limit process becoming mixed Poisson (see, e. g.,
\cite{Grandell1997}). This conclusion corresponds to corollary 3.

\subsection{Generalized Kovalenko theorem. The case of heavy tails}

In this section we will consider the case where condition (22) does
not hold, that is, the tails of the common distribution of the
summands $X_1,X_2,...$ are so heavy that the mathematical
expectation does not exist. Instead of (22) here we will first
assume that ${\sf P}(X_1\ge0)=1$ and $\mathcal{L}(X_1)\in
DNA\big(F(x;\alpha,1)\big)$ for some $\alpha\in(0,1)$, that is,
$$
n^{-1/\alpha}\Sigma_n\Longrightarrow S_{\alpha,1}\ \ \ \
(n\to\infty).\eqno(27)
$$

\smallskip

{\sc Theorem 6.} {\it Assume that i.i.d. nonnegative r.v.'s
$X_1,X_2,...$ satisfy condition $(27)$ with some $\alpha\in(0,1)$.
Let for each $n\in\mathbb{N}$ the r.v. $N_{r,\alpha',1/n^{\alpha'}}$
have the GNB distribution with parameters $r>0$,
$\alpha'\in\mathbb{R}$, $\mu=1/n^{\alpha'}$ and be independent of
the sequence $X_1,X_2,...$ Then
$$
\lim_{n\to\infty}\sup_{x\ge0}\big|{\sf
P}(\Sigma_{N_{r,\alpha',1/n^{\alpha'}}}<n^{1/\alpha}x)-A(x;r,\alpha,\alpha')\big|=0,\eqno(28)
$$
where
$$
A(x;r,\alpha,\alpha')={\sf
P}\big(S_{\alpha,1}G_{r,1}^{1/\alpha\alpha'}<x\big)=\frac{1}{\Gamma(r)}\int_{0}^{x}F(xz^{-1/\alpha\alpha'}\!;\alpha,1)z^{r-1}e^{-z}dz,\
\ \ \ x\ge0. \eqno(29)
$$
If, moreover, $r\in(0,1)$ and $\alpha'\in(0,1]$, then the limit
gamma-mixed strictly stable distribution $A(x;r,\alpha,\alpha')$ is
mixed exponential and hence, infinitely divisible. Moreover, the
corresponding limit r.v. can be represented as
$$
S_{\alpha,1}G_{r,1}^{1/\alpha\alpha'}\eqd
S_{\alpha,1}(G^*_{r,\alpha',1})^{1/\alpha}\eqd
\frac{W_1^{1/\alpha}S_{\alpha,1}}{S_{\alpha',1}^{1/\alpha}Z_{r,1}^{1/\alpha\alpha'}}\eqd
\frac{W_{\alpha}S_{\alpha,1}}{S_{\alpha',1}^{1/\alpha}Z_{r,1}^{1/\alpha\alpha'}}\eqd
\frac{R_{\alpha}W_1}{S_{\alpha',1}^{1/\alpha}Z_{r,1}^{1/\alpha\alpha'}}\eqd
\frac{M_{\alpha}}{S_{\alpha',1}^{1/\alpha}Z_{r,1}^{1/\alpha\alpha'}}\eqd
$$
$$
\eqd
\frac{|X|\sqrt{2W_1}R_{\alpha}}{S_{\alpha',1}^{1/\alpha}Z_{r,1}^{1/\alpha\alpha'}}\eqd
\frac{|\Lambda|R_{\alpha}}{S_{\alpha',1}^{1/\alpha}Z_{r,1}^{1/\alpha\alpha'}}\eqd
S_{\alpha,1}\bigg(\frac{W_{\alpha'}}{Z_{r,1}^{1/\alpha'}}\bigg)^{1/\alpha}\eqd
\frac{S_{\alpha,1}W_{\alpha\alpha'}}{Z_{r,1}^{1/\alpha\alpha'}},
\eqno(30)
$$
where in all terms the factors are independent.}

\smallskip

{\sc Proof.} From condition (27) it follows that in lemma 5 we can
take $b_n\equiv n^{1/\alpha}$, $Q_n\equiv \Sigma_n$, $Q\equiv
S_{\alpha,1}$. Theorem 3 implies relation (25) so that we can take
$d_n\equiv n^{1/\alpha}$, $V\eqd (G^*_{r,\alpha',1})^{1/\alpha}$. By
virtue of lemma 5 conditions (27) and (25) imply
$$
n^{-1/\alpha}\Sigma_{N_{r,\alpha',1/n^{\alpha'}}}\Longrightarrow
S_{\alpha,1}(G^*_{r,\alpha',1})^{1/\alpha}\eqd
S_{\alpha,1}G_{r,1}^{1/\alpha\alpha'}\eqno(31)
$$
as $n\to\infty$. The observation that the limit scale mixture of a
strictly stable distribution is absolutely continuous leads to the
conclusion that the convergence in distribution (31) is uniform over
$x\ge0$. Relation (28) is thus proved.

In the case $r\in(0,1)$ and $\alpha\in(0,1)$ representations (30)
can be easily obtained from theorem 1, representation (6), lemma 3
and the easily verified relation $W_1\eqd
|X|\sqrt{2W_1}\eqd|\Lambda|$. In this case the infinite divisibility
of the limit law follows from the representability of the latter as
a mixed exponential distribution (the fifth term of (30)) by virtue
of a result on generalized gamma-convolutions in
\cite{Bondesson1979} (see corollary 2 there) or a result of Goldie
\cite{Goldie1967} stating that the product of two independent
non-negative random variables is infinitely divisible if one of the
two is exponentially distributed. The theorem is proved.

\smallskip

{\sc Remark 6.} According to theorem 2, if $r\in(0,1)$ and
$\alpha'\in(0,1)$, then the distribution of
$N_{r,\alpha',1/n^{\alpha'}}$ is mixed geometric and the convergence
$$
n^{-1/\alpha}\Sigma_{N_{r,\alpha',1/n^{\alpha'}}}\Longrightarrow
\frac{S_{\alpha,1}W_{\alpha}}{S_{\alpha',1}^{1/\alpha}Z_{r,1}^{1/\alpha\alpha'}}
$$
can be obtained as a corollary of theorem 3 of \cite{Korolev2016TVP}
establishing the conditions for the convergence of the distributions
of mixed geometric random sums.

\smallskip

{\sc Remark 7.} The limit distribution (29) is heavy-tailed: its
moments of orders $\beta\ge\alpha$ do not exist. However, with the
account of the formula expressing the moments of strictly stable
distributions (see the Introduction), for any $\alpha\in(0,1)$,
$\beta\in(0,\alpha)$ and $\alpha'>-\beta/(\alpha r)$ we have
$$
\int_{0}^{\infty}x^{\beta}dA(x;r,\alpha,\alpha')={\sf
E}S_{\alpha,1}^{\beta}{\sf
E}G_{r,1}^{\beta/\alpha\alpha'}=\frac{\Gamma\big(1-\frac{\beta}{\alpha}\big)\Gamma\big(r+\frac{\beta}{\alpha\alpha'}\big)}{\Gamma(1-\beta)\Gamma(r)}.
$$
This relation can be used for the numerical evaluation of
statistical estimates of the parameters $r$,$\alpha$ and $\alpha'$
by the method of moments, say, in the way it was done in
\cite{Kozubowski2001} with respect to the Mittag-Leffler and Linnik
distributions under some additional conditions.

\smallskip

Now consider the case where the summands $X_1,X_2,...$ can take
values of both signs and instead of (27) here we will assume that
$\mathcal{L}(X_1)\in DNA\big(F(x;\alpha,0)\big)$ for some
$\alpha\in(0,2)$, that is,
$$
n^{-1/\alpha}\Sigma_n\Longrightarrow S_{\alpha,0}\ \ \ \
(n\to\infty).\eqno(32)
$$

\smallskip

{\sc Theorem 7.} {\it Assume that i.i.d. r.v.'s $X_1,X_2,...$
satisfy condition $(32)$ with some $\alpha\in(0,2)$. Let for each
$n\in\mathbb{N}$ the r.v. $N_{r,\alpha',1/n^{\alpha'}}$ have the GNB
distribution with parameters $r>0$, $\alpha'\in\mathbb{R}$,
$\mu=1/n^{\alpha'}$ and be independent of the sequence $X_1,X_2,...$
Then
$$
\lim_{n\to\infty}\sup_{x\ge0}\big|{\sf
P}(\Sigma_{N_{r,\alpha',1/n^{\alpha'}}}<n^{1/\alpha}x)-H(x;r,\alpha,\alpha')\big|=0,\eqno(33)
$$
where
$$
H(x;r,\alpha,\alpha')={\sf
P}\big(S_{\alpha,0}G_{r,1}^{1/\alpha\alpha'}<x\big)=
\int_{0}^{\infty}\Phi\Big(\frac{x}{\sqrt{2z}}\Big)dA(z;\alpha/2,\alpha'),\
\ \ x\in\mathbb{R},\eqno(34)
$$
and the d.f. $A(z;r,\alpha/2,\alpha')$ is defined in $(29)$.}

{\it If $r\in(0,1)$ and $\alpha'\in(0,1]$, then the limit normal
scale mixture $H(x;r,\alpha,\alpha')$ is infinitely divisible and
the corresponding limit r.v. can be represented as the scale mixture
of Laplace distributions$:$
$$
S_{\alpha,0}G_{r,1}^{1/\alpha\alpha'}\eqd\frac{\Lambda}{Z_{r,1}^{1/\alpha\alpha'}}\cdot\sqrt{\frac{S_{\alpha/2,1}}{S_{\alpha\alpha'/2,1}}}.
$$
If, in addition, $\alpha'=1$ $($that is, if the r.v. $N_{r,1,1/n}$
has the negative binomial distribution$)$, then the limit r.v. can
be represented as
$$
S_{\alpha,0}G_{r,1}^{1/\alpha}\eqd\frac{\Lambda\sqrt{R_{\alpha/2}}}{Z_{r,1}^{1/\alpha}}\eqd
\frac{X\sqrt{2M_{\alpha/2}}}{Z_{r,1}^{1/\alpha}}\eqd\frac{L_{\alpha}}{Z_{r,1}^{1/\alpha}}.\eqno(35)
$$
where in all terms the factors are independent.}

\smallskip

{\sc Proof.} By the same reasoning that was used to prove theorem 6
we conclude that
$$
\lim_{n\to\infty}\sup_{x\ge0}\bigg|{\sf
P}(\Sigma_{N_{r,\alpha',1/n^{\alpha'}}}<n^{1/\alpha}x)-
\frac{1}{\Gamma(r)}\int_{0}^{x}F(xz^{-1/\alpha\alpha'};\alpha,0)z^{r-1}e^{-z}dz\bigg|=0.
$$
Here the limit d.f. corresponds to the r.v.
$S_{\alpha,0}G_{r,1}^{1/\alpha\alpha'}$. Using relation (1) we
obtain
$$
S_{\alpha,0}G_{r,1}^{1/\alpha\alpha'}\eqd
X\sqrt{2S_{\alpha/2,1}G_{r,1}^{2/\alpha\alpha'}}.\eqno(36)
$$
Notice that the r.v. under the square root sign on the right-hand
side of (36) coincides with the limit r.v. in theorem 6 with
$\alpha$ replaced by $\alpha/2$ (see (29)), whence we conclude that
(33) and (34) hold. In theorem 6 we proved that if $r\in(0,1)$ and
$\alpha'\in(0,1]$, then the d.f. $A(x;r,\alpha,\alpha')$ is
infinitely divisible. It is well known that a normal scale mixture
is infinitely divisible, if the mixing distribution is infinitely
divisible (see, e. g., \cite{Feller}). Hence, if $r\in(0,1)$ and
$\alpha'\in(0,1]$, then the normal scale mixture
$H(x;r,\alpha,\alpha')$ is infinitely divisible.

In the case $r\in(0,1)$ and $\alpha'\in(0,1]$ representations (35)
can easily be obtained from (36), lemma 1 and representations (6)
and (9). The theorem is proved.

\smallskip

{\sc Remark 8.} In the formulation of theorem 7 the representation
of the limit law as a normal scale mixture was chosen just for
reasons of simplicity and convenience. If, instead of (1), the
general form of the <<multiplication theorem>> for strictly stable
laws (see the Introduction) is used, then instead of (36), for any
$\alpha^*\in(\alpha,2]$ we obtain the representation
$$
S_{\alpha,0}G_{r,1}^{1/\alpha\alpha'}\eqd
S_{\alpha^*\!,\,0}\big(S_{\alpha/\alpha^*\!,1}G_{r,1}^{\alpha^*\!/\alpha\alpha'}\big)^{1/\alpha^*}\!,
$$
so that the limit law can be represented as
$$
H(x;r,\alpha,\alpha')={\sf
P}\big(S_{\alpha,0}G_{r,1}^{1/\alpha\alpha'}<x\big)=
\int_{0}^{\infty}F\big(xz^{-1/\alpha^*}\!;\alpha^*\!,0\big)dA(z;r,\alpha/\alpha^*\!,\,\alpha'),\
\ \ x\in\mathbb{R}.
$$

\smallskip

{\sc Remark 9.} According to theorem 2, if $r\in(0,1)$ and
$\alpha'\in(0,1)$, then the distribution of
$N_{r,\alpha',1/n^{\alpha'}}$ is mixed geometric and the convergence
$$
n^{-1/\alpha}\Sigma_{N_{r,\alpha',1/n^{\alpha'}}}\Longrightarrow
S_{\alpha,0}G_{r,1}^{1/\alpha\alpha'}\eqd
\frac{S_{\alpha,0}W_{\alpha}}{S_{\alpha',1}^{1/\alpha}Z_{r,1}^{1/\alpha\alpha'}}\eqd
\frac{S_{\alpha,0}W_{\alpha\alpha'}}{Z_{r,1}^{1/\alpha\alpha'}}
$$
can be obtained as a corollary of theorem 4 of \cite{Korolev2016TVP}
establishing the conditions for the convergence of the distributions
of mixed geometric random sums.

\smallskip

{\sc Remark 10.} The limit distribution (34) is heavy-tailed: its
moments of orders $\beta\ge\alpha$ do not exist. However, with the
account of the formula expressing the moments of strictly stable
distributions (see the Introduction), for any $\alpha\in(0,2)$,
$\beta\in(0,\alpha)$ and $\alpha'>-\beta/(\alpha r)$ we have
$$
\int_{-\infty}^{\infty}|x|^{\beta}dH(x;r,\alpha,\alpha')={\sf
E}|S_{\alpha,0}|^{\beta}{\sf
E}G_{r,1}^{\beta/\alpha\alpha'}=\frac{2^{\beta}}{\sqrt{\pi}}\cdot
\frac{\Gamma\big(\frac{\beta+1}{2}\big)\Gamma\big(\frac{\alpha-\beta}{2}\big)
\Gamma\big(r+\frac{\beta}{\alpha\alpha'}\big)}{\Gamma\big(\frac{2-\beta}{\beta}\big)\Gamma(r)}.
$$
This relation can be used for the numerical evaluation of the
statistical estimates of the parameters $r$,$\alpha$ and $\alpha'$
by the method of moments.

\subsection{An analog of the central limit theorem for GNB random sums}

Consider a sequence of independent identically distributed (i.i.d.)
r.v.'s $X_1,X_2,...$ Assume that ${\sf E}X_1=0$, $0<\sigma^2={\sf
D}X_1<\infty$. For a natural $n\ge1$ let ${S}_n=X_1+...+X_n$. Let
$N_1,N_2,...$ be a sequence of nonnegative integer random variables
defined on the same probability space so that for each $n\ge1$ the
random variable $N_n$ is independent of the sequence $X_1,X_2,...$ A
random sequence $N_1,N_2,...$ is said to be infinitely increasing
($N_n\longrightarrow\infty$) in probability, if ${\sf P}(N_n\le
m)\longrightarrow 0$ as $n\to\infty$ for any $m\in(0,\infty)$.

\smallskip

{\sc Lemma 6}. {\it Assume that r.v.'s $X_1,X_2,...$ and
$N_1,N_2,...$ satisfy the conditions specified above and
$N_n\longrightarrow\infty$ in probability as $n\to\infty$. A d.f.
$F(x)$ such that
$$
{\sf P}\big({S}_{N_n}<x\sigma\sqrt{n}\big) \Longrightarrow F(x)\ \ \
(n\to\infty),
$$
exists if and only if there exists a d.f. $Q(x)$ satisfying the
conditions} $Q(0)=0$,
$$
F(x)=\int_{0}^{\infty}\Phi\big(x/\sqrt{y}\big)dQ(y),\ \
x\in\mathbb{R},\ \ \ {\sf P}(N_n<nx)\Longrightarrow Q(x) \ \
(n\to\infty).
$$\smallskip

{\sc Proof}. This lemma is a particular case of a result proved in
\cite{Korolev1994}, the proof of which is, in turn, based on general
theorems on convergence of superpositions of independent random
sequences \cite{Korolev1996}. Also see \cite{GnedenkoKorolev1996},
theorem 3.3.2.

\smallskip

Re-denote $n=\mu^{-1/\alpha}$. Then $\mu=1/n^{\alpha}$. Consider the
r.v. $N_{r,\alpha,1/n^{\alpha}}$. From theorem 3 it follows that
$N_{r,\alpha,1/n^{\alpha}}\to\infty$ in probability and
$$
\frac{N_{r,\alpha,1/n^{\alpha}}}{n}\Longrightarrow
\frac{W_1}{S_{\alpha,1}Z_{r,1}^{1/\alpha}}\eqd
\bigg(\frac{W_1G_{r,1}}{G_{r,1}+G_{1-r,1}}\bigg)^{1/\alpha}\eqd
W_{\alpha}\cdot\big(1+{\textstyle\frac{1-r}{r}Q_{1-r,r}}\big)^{-1/\alpha}\eqno(37)
$$
as $n\to\infty$, where in each term the involved r.v.'s are
independent.

Now from (37), lemma 6 with $N_n=N_{r,\alpha,1/n^{\alpha}}$, theorem
1 and the well-known relation $\Lambda\eqd X\sqrt{2W_1}$ we directly
obtain

\smallskip

{\sc Theorem 8.} {\it Assume that random variables $X_1,X_2,...$ and
$N_{r,\alpha,1/n^{\alpha}}$, $n\in\mathbb{N}$, satisfy the
conditions specified above.

\noindent $(i)$ Let $r>0$, $\alpha\in\mathbb{R}$. Then
$$
\frac{{S}_{N_{r,\alpha,1/n^{\alpha}}}}{\sigma\sqrt{n}}\Longrightarrow
X\cdot\sqrt{G^*_{r,\alpha,\mu}},
$$
as $n\to\infty$, where the r.v.'s $X$ and $G^*_{r,\alpha,\mu}$ are
independent.

\noindent$(ii)$ If, in addition, $r\in(0,1]$, $\alpha\in(0,1]$, then
the limit law can be represented as
$$
X\cdot\sqrt{G^*_{r,\alpha,\mu}}\eqd
X\cdot\sqrt{\frac{W_1}{S_{\alpha,1}Z_{r,\mu}^{1/\alpha}}}\eqd
\frac{\Lambda}{\mu^{1/2\alpha}\sqrt{2S_{\alpha,1}Z_{r,1}^{1/\alpha}}},
$$
where in each term the involved r.v.'s are independent.}

\smallskip

{\sc Remark 11.} If $\beta\ge 0$ and $\beta/\alpha+r>0$, then the
moment of order $\beta\in(0,\infty)$ of the limit distribution
exists and
$$
{\sf E}\big(|X|\sqrt{G^*_{r,\alpha,\mu}}\big)^{\beta}=
\frac{2^{\beta/2}\Gamma\big(\frac{\beta+1}{2}\big)\Gamma\big(\frac{\beta}{\alpha}+r\big)}{\sqrt{\pi}\mu^{\beta/\alpha}\Gamma(r)}.
$$

\smallskip

{\sc Remark 12.} If $r=\alpha=1$, then the limit law in Theorem 8 is
the ``pure'' Laplace distribution of the r.v. $\mu^{-1/2}\Lambda$.

\section{Limit theorems for statistics constructed from samples with random
sizes having the GNB distributions}

In classical problems of mathematical statistics, the size of the
available sample, i. e., the number of available observations, is
traditionally assumed to be deterministic. In the asymptotic
settings it plays the role of infinitely increasing {\it known}
parameter. At the same time, in practice very often the data to be
analyzed is collected or registered during a certain period of time
and the flow of informative events each of which brings a next
observation forms a random point process. Therefore, the number of
available observations is unknown till the end of the process of
their registration and also must be treated as a (random) {\it
observation}. For example, this is so in insurance statistics where
during different accounting periods different numbers of insurance
events (insurance claims and/or insurance contracts) occur and in
high-frequency financial statistics where the number of events in a
limit order book during a time unit essentially depends on the
intensity of order flows. Moreover, contemporary statistical
procedures of insurance and financial mathematics do take this
circumstance into consideration as one of possible ways of dealing
with heavy tails. However, in other fields such as medical
statistics or quality control this approach has not become
conventional yet although the number of patients with a certain
disease varies from month to month due to seasonal factors or from
year to year due to some epidemic reasons and the number of failed
items varies from lot to lot. In these cases the number of available
observations as well as the observations themselves are unknown
beforehand and should be treated as random to avoid underestimation
of risks or error probabilities.

Therefore it is quite reasonable to study the asymptotic behavior of
general statistics constructed from samples with random sizes for
the purpose of construction of suitable and reasonable asymptotic
approximations. As this is so, to obtain non-trivial asymptotic
distributions in limit theorems of probability theory and
mathematical statistics, an appropriate centering and normalization
of random variables and vectors under consideration must be used. It
should be especially noted that to obtain reasonable approximation
to the distribution of the basic statistics, both centering and
normalizing values should be non-random. Otherwise the approximate
distribution becomes random itself and, for example, the problem of
evaluation of quantiles or significance levels becomes senseless.

In asymptotic settings, statistics constructed from samples with
random sizes are special cases of random sequences with random
indices. The randomness of indices usually leads to that the limit
distributions for the corresponding random sequences are
heavy-tailed even in the situations where the distributions of
non-randomly indexed random sequences are asymptotically normal,
see, e. g., \cite{BeningKorolev2002, BeningKorolev2005,
GnedenkoKorolev1996}.

In this section we will consider the transformation of the limit
laws, if the sample size $n$ is replaced with a random variable.
Consider a sequence of i.i.d. r.v.'s $X_1,X_2,...$ Let $N_1,N_2,...$
be a sequence of nonnegative integer random variables defined on the
same probability space so that for each $n\ge1$ the random variable
$N_n$ is independent of the sequence $X_1,X_2,...$

For $n\ge1$ let $U_n=U_n(X_1,...,X_n)$ be a statistic, that is, a
measurable function of the random variables $X_1,...,X_n$. For each
$n\ge1$ define the random variable $U_{N_n}$ by letting
$U_{N_n}(\omega)=
U_{N_n(\omega)}\left(X_1(\omega),...,X_{N_n(\omega)}(\omega)\right)$
for every elementary outcome $\omega\in\Omega$. We will say that the
statistic $U_n$ is asymptotically normal, if there exists
$\vartheta\in\r$ such that
$$
{\sf
P}\left(\sqrt{n}\bigl(U_n-\vartheta\bigr)<x\right)\Longrightarrow\Phi(x)
\ \ \ (n\to\infty).\eqno{(38)}
$$

\smallskip

{\sc Lemma 7.} {\it Assume that $N_n\longrightarrow\infty$ in
probability as $n\to\infty$. Let the statistic $U_n$ be
asymptotically normal in the sense of $(7)$. Then a distribution
function $F(x)$ such that
$$
{\sf P}\left(\sqrt{n}\bigl(U_{N_n}-\vartheta\bigr)<x\right)
\Longrightarrow F(x)\ \ \ (n\to\infty),
$$
exists if and only if there exists a distribution function $Q(x)$
satisfying the conditions $Q(0)=0$,}
$$
F(x)=\int_{0}^{\infty}\Phi\big(x\sqrt{y}\big)dQ(y),\ \
x\in\mathbb{R},\ \ \ {\sf P}(N_n<nx)\Longrightarrow Q(x) \ \
(n\to\infty).
$$

\smallskip

This lemma is a particular case of theorem 3 in \cite{Korolev1995},
the proof of which is, in turn, based on general theorems on
convergence of superpositions of independent random sequences
\cite{Korolev1996}. Also see \cite{GnedenkoKorolev1996}, theorem
3.3.2.

\smallskip

From (37), lemma 7 with $N_n=N_{r,\alpha,1/n^{\alpha}}$ and theorem
1 with the account of the easily verified property of GG
distributions $(G^*_{r,\alpha,\mu})^{-1}\eqd G^*_{r,-\alpha,\mu}$ we
directly obtain

\smallskip

{\sc Theorem 9}. {\it $(i)$ Let $r>0$, $\alpha\in\mathbb{R}$. Assume
that the statistic $U_n$ is asymptotically normal in the sense of
$(38)$. Then
$$
\sqrt{n}\bigl(U_{N_{r,\alpha,1/n^{\alpha}}}-\vartheta\bigr)\Longrightarrow
 X\cdot\sqrt{G^*_{r,-\alpha,\mu}}
$$
as $n\to\infty$, where the r.v.'s $X$ and $G^*_{r,-\alpha,\mu}$ are
independent.

\noindent$(ii)$ If, in addition, $r\in(0,1]$, $\alpha\in(0,1]$, then
the limit law can be represented as
$$
X\cdot\sqrt{G^*_{r,-\alpha,\mu}}\eqd
X\cdot\sqrt{\frac{S_{\alpha,1}Z_{r,1}^{1/\alpha}}{W_1}}\eqd
S_{2\alpha,0}\cdot\sqrt{\frac{Z_{r,1}^{1/\alpha}}{2W_1}},
$$
where in each term the involved r.v.'s are independent.}

\smallskip

{\sc Remark 13}. The distribution of the limit r.v. in theorem 9 is
a special case of the so-called {\it generalized variance gamma
distributions}, see \cite{KorolevZaks2013}. If $\alpha=1$, then
$S_{\alpha,1}\equiv 1$ and according to lemma 1 the limit law in
theorem 4 turns into that of the r.v. $X\sqrt{Z_{r,1}W_1^{-1}}\eqd
XG_{r,1}^{-1/2}$, that is, the Student distribution with $2r$
degrees of freedom (see \cite{BeningKorolev2005,
KorolevBeningShorgin2011}).

\smallskip

{\sc Remark 14}. It is worth noting that the mixing GG distributions
in the limit normal scale mixtures in theorems 8 and 9 differ only
by the sign of the parameter $\alpha$.

\renewcommand{\refname}{References}

\end{document}